
\magnification = \magstep0
\overfullrule = 2pt
\vsize = 525dd
\hsize = 27.0cc
\topskip = 13dd
\hoffset = 1.8cm
\voffset = 1.8cm
\parindent = 0pt
\parskip = 1ex plus 3pt
\newdimen\Footnoteskip \Footnoteskip = 0pt
\tolerance=200
\emergencystretch=10pt

\message{more fonts: petit and fraktur;}

\font\bfbig = cmbx10 scaled \magstep2   

\font\eightrm = cmr10 scaled 800        
\font\sixrm = cmr7 scaled 850
\font\fiverm = cmr5
\font\eighti = cmmi10 scaled 800
\font\sixi = cmmi7 scaled 850
\font\fivei = cmmi5
\font\eightit = cmti10 scaled 800
\font\eightsy = cmsy10 scaled 800
\font\sixsy = cmsy7 scaled 850
\font\fivesy = cmsy5
\font\eightsl = cmsl10 scaled 800
\font\eighttt = cmtt10 scaled 800
\font\eightbf = cmbx10 scaled 800
\font\sixbf = cmbx7 scaled 850
\font\fivebf = cmbx5

\font\fivefk = eufm5                    
\font\sixfk = eufm7 scaled 850
\font\sevenfk = eufm7
\font\eightfk = eufm10 scaled 800
\font\tenfk = eufm10

\newfam\fkfam
        \textfont\fkfam=\tenfk \scriptfont\fkfam=\sevenfk
                \scriptscriptfont\fkfam=\fivefk

\font\fivebbm  = bbmsl10 scaled 500      
\font\sixbbm   = bbmsl10 scaled 600
\font\sevenbbm = bbmsl10 scaled 700
\font\eightbbm = bbmsl10 scaled 800
\font\tenbbm   = bbmsl10

\newfam\bbmfam
        \textfont\bbmfam=\tenbbm \scriptfont\bbmfam=\sevenbbm
                \scriptscriptfont\bbmfam=\fivebbm
\def\bbm{\fam\bbmfam}

\def\eightpoint{%
        \textfont0=\eightrm \scriptfont0=\sixrm \scriptscriptfont0=\fiverm
                \def\rm{\fam0\eightrm}%
        \textfont1=\eighti  \scriptfont1=\sixi  \scriptscriptfont1=\fivei
                \def\oldstyle{\fam1\eighti}%
        \textfont2=\eightsy \scriptfont2=\sixsy \scriptscriptfont2=\fivesy
        \textfont\itfam=\eightit \def\it{\fam\itfam\eightit}%
        \textfont\slfam=\eightsl \def\sl{\fam\slfam\eightsl}%
        \textfont\ttfam=\eighttt \def\tt{\fam\ttfam\eighttt}%
        \textfont\bffam=\eightbf \scriptfont\bffam=\sixbf
                \scriptscriptfont\bffam=\fivebf \def\bf{\fam\bffam\eightbf}%
        \textfont\fkfam=\eightfk \scriptfont\fkfam=\sixfk
                \scriptscriptfont\fkfam=\fivefk
        \textfont\bbmfam=\eightbbm \scriptfont\bbmfam=\sixbbm
                \scriptscriptfont\bbmfam=\fivebbm
        \rm}
\skewchar\eighti='177\skewchar\sixi='177\skewchar\eightsy='60\skewchar\sixsy='60

\def\fk{\fam\fkfam}
\def\petit{\eightpoint}%


\def\today{\ifcase\month\or
        January\or February\or March\or April\or May\or June\or
        July\or August\or September\or October\or November\or December\fi
        \space\number\day, \number\year}
\def\newline{\hfil\break}

\newtoks\date \date={\today} \def\Date#1{\date={#1}}

\long\def\Comment#1\EndComment{${\cal COMMENT:}${\petit#1}}
\long\def\Comment#1\EndComment{}


\newdimen\defaultparindent
\def\ItemA[#1]#2{\defaultparindent=\parindent\parindent=#1%
        \item{#2}\parindent=\defaultparindent\ignorespaces}%
\def\ItemB#1{\defaultparindent=\parindent\parindent=20pt%
        \item{#1}\parindent=\defaultparindent\ignorespaces}%
\def\ItemOPT{\ifx\ITEMX[\let\ItemCMD\ItemA\else\let\ItemCMD\ItemB\fi\ItemCMD}%
\def\Item{\futurelet\ITEMX\ItemOPT}%


\def\OneColumn{\output={%
        \shipout\vbox{\makeheadline%
                \ifnum\pageno=1\hbox to\fullhsize{\box\TITLEBOX\hfil}\fi%
        \hbox to\fullhsize{\pagebody}%
        \makefootline}%
                \ifnum\pageno=1\global\advance\vsize by\TITLEHEIGHT\fi%
        \advancepageno%
        \ifnum\outputpenalty>-20000%
        \else\dosupereject%
        \fi}%
        }%
\OneColumn

\newbox\LEFTCOLUMN
\def\TwoColumn{\let\lr=L \hsize=.485\fullhsize \output={%
        \if L\lr\global\setbox\LEFTCOLUMN=\leftline{\pagebody}\global\let\lr=R%
        \else \shipout\vbox{\makeheadline%
                \ifnum\pageno=1\hbox to\fullhsize{\box\TITLEBOX\hfil}\fi%
        \hbox to\fullhsize{\box\LEFTCOLUMN\hfil\leftline{\pagebody}}%
        \makefootline}%
                \ifnum\pageno=1\global\advance\vsize by\TITLEHEIGHT\fi%
        \advancepageno\global\let\lr=L%
        \fi%
        \ifnum\outputpenalty>-20000%
        \else\dosupereject%
        \fi}%
        }%

\message{ lay out: title page,}

\def\Date#1{}

\newdimen\fullhsize \fullhsize=\hsize

\newbox\TITLEBOX \setbox\TITLEBOX=\vbox{}
\newdimen\TITLEHEIGHT
\long\def\Title#1#2#3#4{\setbox\TITLEBOX=\vbox{%
        \topglue3truecm%
        \noindent{\bfbig#1}\vskip12pt%
        \noindent{\bf#2}\vskip6pt%
        \noindent{\petit#3}\vskip10pt%
        \noindent{\petit\today}\vskip32pt%
        \noindent{\bf Summary.}\enspace#4\vskip32pt
        \relax}%
        \TITLEHEIGHT=\ht\TITLEBOX%
        \if\IFTHANKS Y\advance\TITLEHEIGHT by\ht\THANKBOX\fi%
        \advance\vsize by-\TITLEHEIGHT%
        \relax}

\newtoks\RunAuthor\RunAuthor={} \newtoks\RunTitle\RunTitle={}

\headline={\ifnum\pageno=1 {\hfil} \else
        \ifodd\pageno {\petit{\the\RunTitle}\hfil\folio}
        \else {\petit\folio\hfil{\the\RunAuthor}} \fi \fi}
\def\makeheadline{\vbox{%
        \hbox to\fullhsize{\the\headline}%
        \vss\nointerlineskip\kern2pt%
        \hbox to\fullhsize{\hrulefill}\kern7pt}}

\newbox\THANKBOX \let\IFTHANKS N
\setbox\THANKBOX=\vbox{\kern12pt\hbox to\fullhsize{\hrulefill}\kern4pt}
\def\Thanks#1#2{\nobreak${}^{#1}$\global\let\IFTHANKS Y%
        \global\setbox\THANKBOX=\vbox{\parindent20pt\baselineskip9pt%
        \unvbox\THANKBOX{\petit\item{${}^{#1}$}{#2}} }}
\footline={\ifnum\pageno=1 \if\IFTHANKS Y\box\THANKBOX\fi
        \else\hfill\fi}
\def\makefootline{\hbox to\fullhsize{\the\footline}}

\message{headings and}

\newcount\SECNO \SECNO=0
\newcount\SSECNO \SSECNO=0
\newcount\SSSECNO \SSSECNO=0
\def\Section#1{\SSECNO=0\SSSECNO=0 \advance\SECNO by 1
  \goodbreak\vskip14pt\noindent{\bf\the\SECNO .~#1}%
  \nobreak\vskip8pt\noindent\kern0pt}
\def\SubSection#1{\SSSECNO=0 \advance\SSECNO by 1
  \goodbreak\vskip14pt\noindent{\it\the\SECNO.\the\SSECNO~#1}%
  \nobreak\vskip8pt\noindent\kern0pt}
\def\SubSubSection#1{\advance\SSSECNO by 1
  \goodbreak\vskip14pt\noindent{\rm\the\SECNO.\the\SSECNO.\the\SSSECNO~#1}%
  \nobreak\vskip8pt\noindent\kern0pt}

\message{definitions,}

\long\def\Def[#1]#2{\medbreak\noindent{\bf#1.\enspace}%
        {\it#2}\medbreak\smallskip\relax}

\long\def\Proof#1\EndProof{\noindent{\it Proof.}\enspace
        #1\hfill$\triangleleft$\medbreak\smallskip\relax}

\def\Verbatim{\bgroup%
        \let\Par=\par\def\par{\Par\leavevmode\relax}%
        \obeylines\obeyspaces\tt%
        \lineskip=0pt\baselineskip=0pt%
        \parindent=10pt\def\\{\tt\char092}%
        \catcode`\%=12 \catcode`\#=12 \catcode`\&=12
        \catcode`\$=12 \catcode`\_=12 \catcode`\^=12
        \catcode`\~=12 \catcode`\{=12 \catcode`\}=12
        \relax}
\def\EndVerbatim{\let\par=\Par\vskip-2\parskip\egroup\relax}


\newcount\FOOTNO \FOOTNO=0
\long\def\Footnote#1{\global\advance\FOOTNO by 1
        {\parindent=16pt\baselineskip=9pt\parskip=-12pt%
        \footnote{\nobreak${}^{\the\FOOTNO)}$}{\petit\hskip-1ex\relax#1\par
        \vskip\Footnoteskip}}}

\message{figures;}

\newdimen\psunit \psunit=1bp
\def\psfile(#1,#2)(#3,#4)[#5,#6]#7{\immediate\write1{}%
        \immediate\write1{psfile data, #7: #5/100 (scale factor)}%
        \immediate\write1{\the\hsize (hsize,text), \the\vsize (vsize,text)}%
        \psunit=#3bp\message{#3bp=\the\psunit (hsize,figure),}%
        \psunit=#4bp\message{#4bp=\the\psunit (vsize,figure)}%
        \vbox to #4bp{\vss\hbox to #3bp{\kern0pt\relax%
        \includegraphics{#7}\hfil}}\ignorespaces}

\newcount\FIGNO \FIGNO=0
\def\Fcaption#1{\global\advance\FIGNO by 1
        {\petit{\bf Fig.~\the\FIGNO.~}#1}}

\def\LArrow[#1]_#2^#3{\mathrel{\smash{\mathop{\hbox to#1{\leftarrowfill}}%
        \limits_{#2}^{#3}}}\relax}
\def\RArrow[#1]_#2^#3{\mathrel{\smash{\mathop{\hbox to#1{\rightarrowfill}}%
        \limits_{#2}^{#3}}}\relax}
\def\UArrow[#1]#2#3{\normallineskip=8pt\matrix{\llap{$\scriptstyle#2$}%
        &\kern-12pt\left\uparrow\vcenter to #1{}\right.&%
        \kern-12pt\rlap{$\scriptstyle#3$}\cr}\relax}
\def\DArrow[#1]#2#3{\normallineskip=8pt\matrix{\llap{$\scriptstyle#2$}%
        &\kern-12pt\left\downarrow\vcenter to #1{}\right.&%
        \kern-12pt\rlap{$\scriptstyle#3$}\cr}\relax}

\message{ referencing system: hlabels,}

\def\htag[#1]{\ignorespaces}
\def\href[#1]#2{#2}

\newtoks\URL \URL={} \def\hlink[#1]#2{%
        \bgroup\edef\~{\string~}\URL={#1}{\href[#1]{#2}}\egroup}

\def\hlabel#1[#2]#3{\htag[#2]\xdef#1{\nobreak\href[#2]{#3}}%
        \ignorespaces}
\def\plabel#1{\hlabel#1[page\the\pageno]{p.\the\pageno}}

\message{equation numbers,}

\newcount\EQNO \EQNO=0
\def\Eqno{\global\advance\EQNO by 1 \eqno(\the\EQNO)%
        \gdef\Label##1{\hlabel##1[eqno\the\EQNO]{(\the\EQNO)}}}

\newcount\THMNO \THMNO=0
\long\def\Thm[#1]#2{\global\advance\THMNO by 1
        \medbreak\noindent{\bf\the\THMNO.~#1.\enspace}%
        {\it#2}\medbreak\smallskip\relax%
        \gdef\Label##1{\hlabel##1[thm\the\THMNO]{\the\THMNO}}}

\message{and the bibliography;}

\newcount\REFNO \REFNO=0
\newbox\REFBOX \setbox\REFBOX=\vbox{}
\def\BegRefs{\message{reading the references}\setbox\REFBOX\vbox\bgroup
        \parindent18pt\baselineskip9pt\petit}
\def\Ref#1{\message{.}\global\advance\REFNO by 1 \ifnum\REFNO>1\vskip3pt\fi
        \item{\the\REFNO .~}\hlabel#1[ref\the\REFNO]{{[\the\REFNO]}}}
\def\EndRefs{\par\egroup\message{done}}
\def\References{\goodbreak\vskip21pt\leftline{\bf References}
        \nobreak\vskip12pt\unvbox\REFBOX\vskip21pt\relax}

\message{some math things.  }

\def\N{I\kern-.8ex N}
\def\Z{\raise.72ex\hbox{${}_{\not}$}\kern-.45ex {\rm Z}}
\def\Q{\raise.82ex\hbox{${}_/$}\kern-1.35ex Q} \def\R{I\kern-.8ex R}
\def\C{\raise.87ex\hbox{${}_/$}\kern-1.35ex C} \def\H{I\kern-.8ex H}
\def\P{I\kern-.8ex P} \def\RP{{\R\!P}}  

\def\N{{\bbm N}}\def\Z{{\bbm Z}}\def\R{{\bbm R}}\def\C{{\bbm C}}
\def\P{{\bbm P}}
\def\RP{\R\!{\sl P}}

\def\D#1#2{{{\partial#1}\over{\partial#2}}}


\Date{\today}
\Title{Harmonic maps in unfashionable geometries}
   {Francis Burstall and Udo Hertrich-Jeromin}
   {Dept.~Math., University of Bath, Bath BA2 7AY,
      email: feb@maths.bath.ac.uk\newline
   Dept.~Math., TU Berlin, Str.~d.~17.~Juni 136, D-10623 Berlin,
      email: udo@sfb288.math.tu-berlin.de}
   {We describe some general constructions on a real smooth projective
   $4$-quadric which provide analogues of the Willmore functional and
   conformal Gauss map in Lie sphere and projective differential
   geometry.  Extrema of these functionals are characterized by
   harmonicity of this Gauss map.}%
   \footnote{}{\petit 2000 Mathematics Subject Classification:
      53C43 (primary), 58E20 53A20 53A40 (secondary)\newline
      }

\BegRefs
\Ref\Blaschke 
   W.~Blaschke: {\it Vorlesungen \"uber Differentialgeometrie III\/};
   Springer Grundlehren 29, Berlin 1929
\Ref\BolII 
   G.~Bol: {\it Projektive Differentialgeometrie, 2.~Teil\/};
   Math.~Lehrbuecher 9, Vandenhoeck \& Rup\-recht, G\"ottingen 1954
\Ref\Bryant
   R.L.~Bryant: {\it {C}onformal and minimal immersions of compact surfaces
   into the {$4$}-sphere\/}; J. Diff.\ Geom.~{\bf 17} (1982) 455-473
\Ref\BFLPP
   F.~Burstall, D.~Ferus, K.~Leschke, F.~Pedit, U.~Pinkall:
   {\it Conformal Geometry of Surfaces is $S^4$ and Quaternions\/};
   Preprint 1999
\Ref\BFPP 
   F.E.~Burstall, D.~Ferus, F.~Pedit, U.~Pinkall:
   {\it {H}armonic tori in symmetric spaces and commuting Hamiltonian
   systems on loop algebras\/}; Ann.\ Math.~{\bf 138} (1993) 173-212
\Ref\BP 
   F.E.~Burstall, F.~Pedit: {\it {D}ressing orbits of harmonic maps\/};
   Duke Math.\ J.~{\bf 80} (1995) 353-382
\Ref\Cecil 
   T.~Cecil: {\it Lie sphere geometry\/}; Springer Universitext,
   New York 1992
\Ref\Demoulin 
   A.~Demoulin: {\it Sur deux transformations des surfaces dont les
   quadriques de Lie n'ont que deux ou trois points caracteristiques\/};
   Bull.\ Acad.\ Roy.\ Belg., V.~Ser.~{\bf19} (1933) 1352-1363
\Ref\FerapontovD 
   E.~Ferapontov, W.~Schief: {\it Surfaces of Demoulin:
   differential geometry, B\"acklund transformation, and integrability\/};
   J.~Geom.\ Phys.~{\bf30} (1999) 343-363
\Ref\FerapontovL
   E.~Ferapontov: {\it Lie sphere geometry and integrable
   systems\/}; Tohoku Math.~J.~{\bf52} (2000) 199-233
\Ref\FerapontovP
   E.~Ferapontov: {\it Integrable systems in projective
   differential geometry\/}; Kyushu J.~Math.~{\bf54} (2000) 183-215
\Ref\Godeaux 
   L.~Godeaux: {\it Sur une suite de surfaces dont les quadriques
   de Lie n'ont que trois points caracteristiques\/};
   Bull.\ Acad.\ Roy.\ Belg., V.~Ser.~{\bf20} (1934) 495-504
\Ref\Helgason 
   S.~Helgason: {\it Differential geometry, Lie groups, and symmetric
   spaces\/}; Academic Press, New York 1978
\Ref\Hitchin 
   N.J.~Hitchin: {\it Harmonic maps from a $2$-torus to the $3$-sphere\/};
   J.\ Diff.\ Geom.~{\bf 31} (1990) 627-710
\Ref\Klein 
   F.~Klein: {\it Vorlesungen \"uber h\"ohere Geometrie\/};
   Springer Grundlehren 22, Berlin 1926
\Ref\LieBrief
   S.~Lie, letter to F.~Klein, 18.12.1878; collected papers, vol.~3
   (1922) 718-719
\Ref\Lie 
   S.~Lie: {\it Geometrie der Ber\"uhrungstransformationen\/};
   Chelsea, New York 1977
\Ref\Liebermann 
   P.~Libermann, C.-M.~Marle: {\it Symplectic geometry and analytical
   mechanics\/}; Mathematics and its Applications 35, Kluwer Academic
   Publishers, Dordrecht 1987
\Ref\Pohlmeyer
   K.~Pohlmeyer: {\it {I}ntegrable {H}amiltonian systems and interactions
   through quadratic constraints\/}; Commun.\ Math.\ Phys.~{\bf 46} (1976)
   207-221
\Ref\Rozet 
   O.~Rozet: {\it Sur certaines congruences $W$ attachees aux surfaces
   dont les quadriques de Lie n'ont que deux points caracteristiques\/};
   Bull.\ Sci.\ Math., II.~Ser.~{\bf58} (1934) 141-151
\Ref\TerngUhlenbeck
   C.-L.~Terng, K.~Uhlenbeck: {\it {B}\"acklund transformations and
   loop group actions\/};  Comm.\ Pure Appl.\ Math.~{\bf 53} (2000) 1-75
\Ref\ThomsenA 
   G.~Thomsen: {\it \"Uber eine liniengeometrische Behandlungsweise
   der projektiven Fl\"achentheorie und die projektive Geometrie der Systeme
   von Fl\"achen zweiter Ordnung\/}; Abh.\ Math.\ Sem.\ Hamburg {\bf4} (1925)
   232-266
\Ref\ThomsenB 
   G.~Thomsen: {\it Sulle superficie minime proiettive\/}; Annali
   Mat.~{\bf5} (1928) 169-184
\Ref\Uhlenbeck
   K.~Uhlenbeck: {\it {H}armonic maps into {L}ie groups (classical solutions
   of the chiral model)\/}; J. Diff.\ Geom.~{\bf 30} (1989) 1-50
\Ref\ZSh
   V.E.~Zakharov, A.B.~Shabat: {\it Integration of non-linear equations of
   mathematical physics by the method of inverse scattering, II\/};
   Funkts.\ Anal.\ i Prilozhen~{\bf 13} (1978) 13-22
\EndRefs

\def\f{{\fk f}}

\Section{Introduction}
Many topics in integrable surface geometry\Footnote{But not all:
isothermic surfaces, for example, do not fit into this picture.} may
be unified by application of the highly developed theory of harmonic
maps of surfaces into (pseudo-)Riemannian symmetric spaces.  On the
one hand, such harmonic maps comprise an integrable system with
spectral deformations, algebro-geometric solutions and dressing
actions of loop groups generated by B\"acklund transforms \BFPP, \BP,
\Hitchin, \TerngUhlenbeck, \Uhlenbeck.  On the other hand, several
integrable classes of surface are characterized by harmonicity of a
suitable Gauss map.  Thus, a surface $\f:M^2\to\R^3$ has constant
mean curvature $H$ if and only if its Gauss map $M^2\to S^2$ is
harmonic.  Again, such a surface has constant Gauss curvature $K$ if
and only if its Gauss map is harmonic with respect to the metric on
$M$ provided by the second fundamental form of $\f$.  The theory
of harmonic maps now provides a conceptual explanation of the classical
integrable aspects of such surfaces such as associated families, Lie
and B\"acklund transformations.

These ideas gain wider applicability if we extend the notion of
Gauss map.  Consider, for example, the case of Willmore surfaces
\Blaschke: these are surfaces $\f:M^2\to\R^3$ which extremize the
Willmore functional
$$
W(\f)=\int_{M^2}(H^2-K)\, dA.
$$
The functional $W$ and so its critical points are preserved by the
M\"obius group of conformal diffeomorphisms of
$S^3=\R^3\cup\{\infty\}$.  A conformally immersed surface in
$S^3$ also has a Gauss map which can be defined as follows: to each
point $x\in M^2$, attach the oriented $2$-sphere $S(x)$ in $S^3$ which has
first order contact with $\f$ at $\f(x)$ and the same mean curvature
vector there.  The map $x\mapsto S(x)$ is variously known as the {\it
central sphere congruence\/} \Blaschke\ or the {\it conformal Gauss
map\/} \Bryant\ and is a M\"obius invariant of $\f$.  The space of
oriented $2$-spheres in $S^3$ is naturally identified with the
Lorentz $4$-sphere which is a pseudo-Riemannian symmetric space.  One
shows that the harmonic map energy $E(S)$ of $S$ coincides with
$W(\f)$ and further \Blaschke, \Bryant\ that $\f$ is Willmore if and
only if $S$ is harmonic.

In this paper we will study two other classes of surfaces, Lie minimal
and projectively minimal surfaces, and show that an exactly analogous
theory applies.  These surfaces are the analogues of Willmore
surfaces in projective and Lie sphere geometry and were introduced
and intensively studied around the turn of the last century (see, for
example, \ThomsenA, \ThomsenB, \Blaschke, \BolII).  More recently,
Ferapontov \FerapontovL, \FerapontovP\ has demonstrated that these
surfaces have integrable structure and we shall show how this
structure is explained by the harmonicity of a Gauss map.

Let us begin by explaining what Lie and projectively minimal surfaces
are.  First contemplate an immersed surface $\f:M^2\to\R^3$ with
curvature line coordinates $u,v$ and corresponding principal
curvatures $\kappa_1,\kappa_2$.  We define a functional $L_{\rm Lie}$
using a formulation we learned from Ferapontov \FerapontovL:
$$
L_{\rm Lie}(\f)=\int_{M^2}{\partial_u\kappa_1\partial_v\kappa_2\over
(\kappa_1-\kappa_2)^2}du\wedge dv.
$$
The critical points of $L_{\rm Lie}$ are called {\it Lie minimal\/}
surfaces.  One can show (and we will!) that the Lagrangian density
(and so $L_{\rm Lie}$ and its critical points) is invariant under
both the M\"obius group and normal shifts (the passage to a parallel
surface ${\fk f}+t{\fk n}$).  Otherwise said, the density is
preserved by the group of Lie sphere transformations (see Cecil
\Cecil\ for a modern account of Lie sphere geometry).

Secondly, let $\f:M^2\to\RP^3$ be an immersed surface in real
projective $3$-space with (possibly complex conjugate) asymptotic
coordinates\Footnote{Note that the notion of asymptotic coordinates
is projectively invariant as the conformal class of the second
fundamental form is.} $u,v$.  Thus we have\Footnote{Note that here,
and elsewhere, we do not distinguish between a map ${\fk
f}:M^2\to\RP^3$ and any lift (expression in homogeneous coordinates)
${\fk f}:M^2\to\R^4_{\times}$.}
$$
\matrix{
   {\fk f}_{uu} &=& *\,{\fk f}_u &+& p\,{\fk f}_v &+& *\,{\fk f}, \cr
   {\fk f}_{vv} &=& q\,{\fk f}_u &+& *\,{\fk f}_v &+& *\,{\fk
f}, \cr}
$$ for some functions $p,q$ (here and elsewhere, we use $*$ to
represent unknown functions that are irrelevant to our analysis).  It
is not difficult to check that the density $pq\,du\wedge dv$ is
independent of choices (both of asymptotic coordinates and lift) so
that we have a well-defined functional
$$
L_{\rm proj}(\f)=\int_{M^2}pq\,du\wedge dv.
$$
The critical points are the {\it projectively minimal surfaces\/}
\ThomsenB.  In this case the density is invariant under the
projective action of ${\rm SL}(4,\R)$ on $\RP^3$.

Our contention is that Lie minimality and projective minimality are
characterized by harmonicity of an appropriate Gauss map and that,
moreover, this Gauss map has a geometric interpretation as a
congruence of ``model surfaces'' --- either Dupin cyclides or
quadrics --- having second order contact with the immersion $\f$.  Both
Dupin cyclides and quadrics of fixed signature are parametrized by
pseudo-Riemannian symmetric spaces (in fact, Grassmannians) and, in
classical language, our main results are:
\Thm[Theorem] {\it $\f:M^2\to\R^3$ is Lie minimal if and only if its
congruence of Lie cyclides is harmonic.\/}
\Thm[Theorem] {\it $\f:M^2\to\RP^3$ is projectively minimal if and only if its
congruence of Lie quadrics is harmonic.\/}

We find a uniform treatment of these assertions in the following
considerations: first we treat the contact lifts of immersions rather
than the immersions themselves --- indeed, for Lie sphere geometry,
this is compulsory since the symmetry group of the situation does not
act by point-transformations: for example, a circle and a torus of
revolution are Lie sphere equivalent via a normal shift.  Second, we
exploit the fact that the space of contact elements in
$S^3=\R^3\cup\{\infty\}$ and in $\RP^3$ share a common description as
the space $Z$ of lines in a $4$-dimensional quadric ${\cal
Q}\subset\RP^5$.  For Lie sphere geometry, this comes from the fact
that oriented $2$-spheres in $S^3$ (including points) are
parametrized by the {\it Lie quadric\/}: the projective light cone of
$\R^{4,2}$ and lines in this quadric correspond to parabolic pencils
of spheres or, equivalently, contact elements in $S^3$ (see \Cecil).
In projective geometry, the double cover ${\rm SL}(4,\R)\to{\rm
O}(3,3)$ gives rise to the {\it Klein correspondence\/} between the
space of lines in $\RP^3$ and the {\it Pl\"ucker quadric\/}: the
projective light cone of $\R^{3,3}$.  Then lines in the Pl\"ucker
quadric parametrize contact elements in $\RP^3$.

Thus we arrive at a uniform approach by considering Legendre
immersions into the space $Z$ of lines in a $4$-dimensional quadric
${\cal Q}$ as described in Section~2.  In this setting, with the aid
of the focal surfaces and conjugate parameters attached to such an
immersion, we shall, in Section~3, equip each Legendre immersion with
a Grassmannian-valued Gauss map.  We shall use this in Section~4 to
define a functional on Legendre immersions whose critical points
(with respect to Legendre variations) are characterized by
harmonicity of this Gauss map.  These constructions proceed
independently of the signature of the metric on $\R^6$ defining
${\cal Q}$ but specialize, as we shall see in Section~5, to give our
main results once that signature is declared.  Thus our methods may
be viewed as a practical implementation of the famous line-sphere
correspondence of Lie\Footnote{Note that Lie's line-sphere
correspondence does {\it not\/} provide an isomorphism between the
Lie and Pl\"ucker quadrics (these two spaces are topologically
different).  However, it can certainly be considered as a
correspondence between concepts in projective line geometry and Lie
sphere geometry (cf.~\Lie, \Klein).}

We conclude our study by indicating in Section~6 some applications of
harmonic map theory to Lie and projectively minimal surfaces.

\Section{Line congruences in quadrics and Legendre surfaces}
We consider a nonsingular 4-dimensional quadric ${\cal Q}\subset\RP^5$
that contains (real) lines, that is, ${\cal Q}=\P({\cal L})$ is the
projectivized light cone of some $\R^{m,n}$ where $m+n=6$ and
$m\geq n\geq2$ (this last inequality is the condition that
${\cal Q}$ contains real lines).
Further, let $Z^{m,n}:=G_{2\times0}(\R^{m,n})$ denote the space of null
2-planes in $\R^{m,n}$, that is, the space of lines in ${\cal Q}$.
In the sequel, we will refer to the space of lines in ${\cal Q}$ as $Z$ unless
the signature of the underlying quadric ${\cal Q}$ has to be emphasized.

The orthogonal group ${\rm O}(m,n)$ acts transitively on these spaces (and
others we shall consider below) and this gives convenient algebraic
models for their tangent spaces which we shall use repeatedly.  This
being the case, let us briefly recall the basic setting of homogeneous
geometry: so let $G$ be a Lie group with Lie algebra ${\fk g}$ and $N$
a homogeneous $G$-space.  Each $\xi\in\fk g$ gives rise to a vector
field $\tilde\xi$ on $N$ via
$\tilde\xi_x=\left. d\over{dt}\right|_{t=0}\exp t\xi\cdot x$
and then
$$ [\tilde\xi,\tilde\eta]=-\widetilde{[\xi,\eta]}.\Eqno\Label\Bracket $$
Since $G$ acts transitively on $N$, we have a surjection
${\fk g}\ni\xi\mapsto\tilde\xi_x\in T_xN$ whose kernel is the (infinitesimal)
stabiliser ${\fk g}^x$ of $x$.
Thus, $T_x N\cong {\fk g}/{\fk g}^x$ and, more globally, we have identified
$TN$ with the quotient of the trivial bundle $N\times\fk g$ by the bundle of
stabilisers.

Let us consider the example of $Z$ in more detail: $\pi\in Z$ is a null
$2$-plane in $\R^{n+m}$ and so ${\fk o}(m,n)^\pi=\{\xi\in{\fk o}(m,n):
\xi\pi\subset\pi\}$.  Thus $T_\pi Z\cong {\fk o}(m,n)/{\fk o}(m,n)^\pi$.
However, it is a simple matter to see that restriction gives a surjection
$$
{\fk o}(m,n)\ni\xi\mapsto\xi_{|\pi}{\rm\ mod}\ \pi\in
\{A\in{\rm Hom}(\pi,\R^{m,n}/\pi)\,|\,
\langle As_1,s_2\rangle+\langle s_1,As_2\rangle=0\}
$$
with kernel ${\fk o}(m,n)^\pi$ so that we have an isomorphism
$$
T_\pi Z\cong
\{A\in{\rm Hom}(\pi,\R^{m,n}/\pi)\,|\,
\langle As_1,s_2\rangle+\langle s_1,As_2\rangle=0\}.
\Eqno\Label\TZ
$$
Explicitly, this isomorphism is given by $X\mapsto A_X$ where
$A_Xs=d_Xs{\rm\ mod}\ \pi$ for any local section $s$ of the tautological
bundle over $Z$ (note that $A_X$ so defined is algebraic).

Two structures on $Z$ will be important to us in the sequel: first $Z$
is a contact manifold.  For this, define a line bundle $L$ by
$L_\pi=\Lambda^2\pi^*$ and note that we have a surjective bundle map
$\vartheta:TZ\to L$ via $\vartheta(A)=\langle
A.,.\rangle|_{\pi\times\pi}$, for $A\in T_\pi Z$.  Denote the kernel
of $\vartheta$ by ${\cal D}$.  Clearly
$$
{\cal D}_\pi={\rm Hom}(\pi,\pi^{\perp}/\pi).
$$ We claim that $\vartheta$ provides a contact structure, that is,
the (algebraic) map ${\cal D}\times{\cal D}\to L$ given by
$X,Y\mapsto\vartheta([X,Y])$ is non-degenerate.  In fact, using
\Bracket, one sees that, for $A,B\in{\cal D}_\pi$,
$$
\vartheta([A,B])(s_1,s_2)=\langle Bs_1,As_2\rangle-\langle
As_1,Bs_2\rangle
\Eqno\Label\Symplectic
$$
which is readily checked to be non-degenerate.

Secondly, both $\pi$ and $\pi^\perp/\pi$ are $2$-dimensional so that
we can equip ${\cal D}$ with a $(2,2)$-conformal structure by setting
$(A,A)=\det A$.  Of course, this requires a choice of bases on $\pi$
and $\pi^\perp/\pi$ but changing this choice merely rescales the
result\Footnote{More invariantly, one can view the conformal structure
as the inner product $A,B\mapsto{1\over2}(A\wedge B+B\wedge A)$ taking
values in the line $\Lambda^2\pi^*\otimes\Lambda^2 \pi^\perp/\pi$.}.

\Def[Definition]{A map $f:M^2\to Z$ is called Legendre if it is
tangent to the contact distribution ${\cal D}$, that is,
if $\langle ds_1,s_2\rangle\equiv0$ for $s_1,s_2:M^2\to{\cal L}$ with
$f=s_1\wedge s_2$.}

A Legendre immersion pulls back the conformal structure on ${\cal D}$ to
one on $M$ that will be useful to us.  Let us consider the possibilities
for its signature: clearly any tangent plane to $f$ is Lagrangian for the
symplectic form \Symplectic\ and, when $(m,n)=(4,2)$, it is easily seen
that this forces the conformal structure on any $T_xM$ to vanish or have
signature $(1,1)$.  When $(m,n)=(3,3)$ there is no such restriction and
all signature are possible.

\Def[Assumption]{From this point on, we assume that $f:M^2\to Z$ is a
Legendre map such that the induced conformal structure ($df$ takes values
in ${\cal D}$!) is non-degenerate.  In particular, wherever applicable, $f$
will be assumed to be an immersion.  Further, let $(u,v)$ be (possibly
complex conjugate) coordinates along the null-directions of the induced
conformal structure.}

As $f_u,f_v\in{\rm Hom}(f,f^{\perp}/f)={\cal D}_f$ are null directions
of the conformal structure on ${\cal D}$ we have $\det f_u=\det f_v=0$
so that there are $l,s:M^2\to{\cal L}$ with ${\rm ker}f_u={\rm
span}\{l\}$ and ${\rm ker}f_v={\rm span}\{s\}$.  Thus,
$l_u,s_v\equiv0$ mod $f$.  In classical language, $l,s$ are the focal
surfaces of the line congruence $f$ with $(u,v)$ forming their common
conjugate net.  

Since $f$ is an immersion, we have $l_v,s_u\not\in f$ as, otherwise,
$f_u=0$ or $f_v=0$.  In particular, $l$ and $s$ are linearly
independent so that $f=l\wedge s$.  Similarly, $l_v$ and $s_u$ are
linearly independent mod $f$ as, otherwise, the induced conformal
structure would have a third null direction.

As we shall see, these constructions have a direct geometric
interpretation: in case $(m,n)=(4,2)$, $f$ is the contact lift of a
surface in $S^3$ for which $(u,v)$ are curvature line
coordinates\Footnote{Curvature lines are real: this explains the
restriction on the signature of the induced conformal structure in
the $(4,2)$-case.} and $l,s$ are the curvature spheres (cf.~\Cecil);
similarly, in case $(m,n)=(3,3)$, $f$ is the contact lift of a
surface in $\RP^3$ for which $(u,v)$ are asymptotic coordinates and
$l,s$ are the corresponding asymptotic line congruences (cf. \Lie).

\Section{The conformal Gauss map}
We first prove that $s,s_u,s_{uu}$ and $l,l_v,l_{vv}$ define two orthogonal
3-dimensional bundles with non-degenerate induced metrics.

Differentiating $\langle l,s\rangle=0$ with respect to $v$ twice
yields $l,l_v,l_{vv}\perp s$ since $s_v\in l\wedge s$; similarly,
$s,s_u,s_{uu}\perp l$.  As $s_{uv}\in l\wedge s\wedge s_u\perp l$, we
find $\langle l_v,s_u\rangle =\langle l,s_u\rangle_v=0$.  Thus, we
also have $s_{uv}\in l\wedge s\wedge s_u\perp l_v$, so that we obtain
$\langle l_{vv},s_u\rangle=\langle l_v,s_u\rangle_v=0$; similarly,
$\langle l_v,s_{uu}\rangle=0$.  Finally, we have $s_{uuv}\in l\wedge
s\wedge s_u\wedge s_{uu}\perp l_v$ which yields $\langle
l_{vv},s_{uu}\rangle=\langle l_v,s_{uu}\rangle_v=0$.  This gives the
first assertion $s,s_u,s_{uu}\perp l,l_v,l_{vv}$.

In particular, $f^{\perp}=l\wedge s\wedge s_u\wedge l_v$ so that the
non-degeneracy of
$\langle,\rangle|_{(f^{\perp}/f)\times(f^{\perp}/f)}$ shows that we
must have $\langle s_u,s_u\rangle,\langle l_v,l_v\rangle\neq0$.  The
pairwise scalar products of $l,l_v,l_{vv}$ and $s,s_u,s_{uu}$,
respectively, are given by $$\matrix{ \left(\matrix{0&0&-\langle
l_v,l_v\rangle\cr 0&\langle l_v,l_v\rangle&\ast\cr -\langle
l_v,l_v\rangle&\ast&\ast\cr}\right) &{\rm and}&
\left(\matrix{0&0&-\langle s_u,s_u\rangle\cr 0&\langle
s_u,s_u\rangle&\ast\cr -\langle s_u,s_u\rangle&\ast&\ast\cr}\right)
\cr}$$ showing that $l,l_v,l_{vv}$ and $s,s_u,s_{uu}$ span
3-dimensional subspaces with non-degenerate metrics.  Note that both
spaces contain null lines and, in case $(m,n)=(4,2)$, these are real
(since $u,v$ are real) so that their signatures are both $(2,1)$.
When $(m,n)=(3,3)$ there are two cases to consider: if $(u,v)$ are
real, both spaces are real and contain null lines so that, as before,
the signatures are $(2,1)$ and $(1,2)$.  If, however, $(u,v)$ are
complex conjugate then these spaces are complex and conjugate to each
other.

In all cases,
we have to do with a splitting $(\R^{m,n})^{\C}=\C^6=S\oplus S^{\perp}$ 
where $S$, $S^{\perp}$ are 3-dimensional orthogonal complex spaces
satisfying a reality condition $$\matrix{
   \bar S = S \hfill&{\rm for}& (u,v) {\rm\ real,\ and}\hfill\cr
   \bar S = S^{\perp} &{\rm for}& (u,v) {\rm\ complex\ conjugate.}\cr}$$
Such a splitting corresponds to a real symmetric endomorphism
$\star_S:\C^6\to\C^6$ with $\star_S^2=\varepsilon^2$, $\varepsilon=1$ or
$=i$, by setting $S$, resp.\ $S^{\perp}$, to be the $+\varepsilon$, resp.\
$-\varepsilon$, eigenspaces of $\star_S$. Now, $\star_S$ gives rise to a
Legendre submanifold $Z_S$ of $Z$ by $$\matrix{
   Z_S := \{\pi\in Z\,|\, \star_S\pi=\pi\}.}$$
Indeed, $T_{\pi}Z_S=\{A\in T_{\pi}Z\,|\,\star_S\!A=A\star_S\}$ and, with
$\varepsilon_{\pm}$-eigenvectors $s_{\pm}$ of $\star_S$ in $\pi^{\C}$, we have
$\langle As_+,s_-\rangle=-{1\over\varepsilon}\langle As_+,\star_Ss_-\rangle
  =-{1\over\varepsilon}\langle A\star_Ss_+,s_-\rangle=-\langle As_+,s_-\rangle$
giving $\langle As_+,s_-\rangle=0$ so that $T_{\pi}Z_S\subset{\cal D}_{\pi}$.
Thus $$\matrix{
  T_{\pi}Z_S=\{A\in{\rm Hom}(\pi,\pi^{\perp}/\pi)\,|\,\star_S\!A=A\star_S\}.}$$
Moreover, the null vectors in $T_{\pi}Z_S$ have $\star_S$-stable kernels and
images and so are real or complex conjugate according to whether $\varepsilon=1$
or $=i$ (the latter case only being possible when $(m,n)=(3,3)$).
Thus $Z_S$ has a conformal structure of signature $(1,1)$ for $\varepsilon=1$
or $(2,0)$ for $\varepsilon=i$.
We label the Grassmannian of all such $S$ by the signature $(m,n)$ of the
real structure and that of the conformal structure on $Z_S$: thus we set
$$\matrix{
   {\cal G}_{i,j}^{m,n} := \{\star_S:\R^{m,n}\to\R^{m,n}\,\left|\,\matrix{
      \star_S{\rm\ is\ symmetric}, \star_S^2=\pm1, {\rm\ and}\cr
      Z_S {\rm\ has\ signature}\ (i,j)}\right\}.}$$
To be absolutely explicit, this means $$\matrix{
   {\cal G}_{1,1}^{m,n} &=& \{S=(S_{\R})^{\C}\subset(\R^{m,n})^{\C}\,|\,
      S_{\R}\subset\R^{m,n}{\rm\ has\ signature}\ (2,1)\}, \hfill\cr
   {\cal G}_{2,0}^{3,3} &=& \{S\subset(\R^{3,3})^{\C}\,|\,
      S\cap S^{\perp}=\{0\}, {\rm\ and}\ \bar S=S^{\perp}\}. \hfill\cr}$$
We will see later that ${\cal G}_{1,1}^{4,2}$ parametrizes contact lifts of
$Z_S$ of Dupin cyclides while ${\cal G}_{i,j}^{3,3}$ parametrizes (contact
lifts of) quadrics in $\RP^3$ with conformal structure of signature $(i,j)$.

In the sequel, we shall drop the decorations and simply refer to any of
these Grassmannians as ${\cal G}$ unless the signatures need emphasis.

\Def[Definition]{Given a Legendre surface $f:M^2\to Z$ with induced conformal
structure of signature $(i,j)$, focal surfaces $l,s:M^2\to{\cal L}^{\C}$, and
conjugate parameters $(u,v)$, we define the conformal Gauss map
$S:M^2\to{\cal G}_{i,j}^{m,n}$ of $f$ by $S:={\rm span}_{\C}\{l,l_v,l_{vv}\}$.}

It is clear that this definition is indeed independent of any choices
(besides swapping the roles of the subspaces in the Lie sphere case).

Here is the geometry of the situation: by construction, $l,s$ are
eigenvectors of $\star_S$ so that $f=l\wedge s$ is $\star_S$-stable.
That is, for each $x\in M$, we have $f(x)\in Z_{S(x)}$.
Additionally, the images $s_u,l_v$ mod $f$ of $f_u,f_v$ are eigenvectors
of $\star_S$ so that $f_u$ and $f_v$ both commute with $\star_S$.
Thus, each $d_xf(T_xM)=T_{f(x)}Z_{S(x)}$ so that $f(M)$ and $Z_{S(x)}$
have first order contact at $f(x)$.  Otherwise said, the congruence
$Z_S$ envelopes $f$.

We now justify our terminology by showing that our conformal Gauss map
$S$ is indeed conformal.  For this, first contemplate the Grassmannians
${\cal G}$: these are pseudo-Riemannian symmetric spaces.
In fact, ${\rm O}(m,n)$ acts isometrically on each and the identification
$T_S{\cal G}$ with ${\fk o}(m,n)/{\fk o}(m,n)^S$ together with the
natural isomorphism $\C^6/S\cong S^{\perp}$ gives an identification
of $T_S^{\C}{\cal G}$ with ${\rm Hom}(S,S^{\perp})$ via
$X\mapsto[\sigma\mapsto\pi^{\perp}(d_X\sigma)]$, where
$\pi^{\perp}:\C^6\to S^{\perp}$ denotes orthogonal projection.
We define a real indefinite\Footnote{It may amuse the reader to
compute the signature.} metric on ${\rm Hom}(S,S^{\perp})$ by
$\langle A,B\rangle=-{\rm tr}B^{\star}A=-{\rm tr}AB^{\star}$ where
$B^{\star}:S^{\perp}\to S$ is the adjoint of $B:S\to S^{\perp}$.
Moreover, $\star_S$ induces an isometric involution on ${\cal G}$
which is the symmetric involution at $S$.

With this in hand, consider a conformal Gauss map $S:M\to{\cal G}$.
We identify $S^{\ast}T^{\C}{\cal G}$ with ${\rm Hom}(S,S^{\perp})$ and
then $dS=[\sigma\mapsto dS\cdot\sigma=\pi^{\perp}d\sigma]$ for $\sigma$
a local section of the bundle $S\to M$.
Note that, for $\sigma\in\Gamma(S)$ and $\varrho\in\Gamma(S^{\perp})$,
we have $\langle dS\cdot\sigma,\varrho\rangle
   =\langle d\sigma,\varrho\rangle
   =-\langle \sigma,d\varrho\rangle$ so that
$dS^{\star}$ is given by $dS^{\star}\varrho=-\pi d\varrho$ where
$\pi$ is the orthogonal projection $\C^6\to S$.

By construction, $S_vl=S_vl_v=0$, that is,
${\rm ker}\,S_v\supset l^{\perp}\cap S$, hence
${\rm im}\,S_v^{\star}\subset{\rm span}\{l\}$.
Similarly, $s^{\perp}\cap S^{\perp}\subset{\rm ker}\,S_u^{\star}$,
or ${\rm im}\,S_u\subset{\rm span}\{s\}$.
Consequently, $S_u^{\star}\circ S_u=0$ and $S_v\circ S_v^{\star}=0$
and taking traces gives $\langle S_u,S_u\rangle=\langle S_v,S_v\rangle=0$.
Since $(u,v)$ are null coordinates for the conformal structure induced
by $f$ we conclude that $S$ is indeed conformal.

\Thm[Theorem]{The conformal Gauss map $S$ of a Legendre surface
$f$ is conformal.}

For later use, we note that $S_v^{\star}\circ S_u(l)=-\langle
S_u,S_v\rangle l$ since ${\rm im}\,S_v^{\star}\subset{\rm
span}\{l\}$.  Thus, defining functions $p,q$ on $M$ by
$$
\matrix{ 
l_u &=& \ast\, l &+& p\,s \cr s_v &=&
q\,l &+& \ast\, s \cr }
$$
yields
$$ \langle S_u,S_v\rangle=pq \Eqno\Label\WillmoreDensity$$ since
$S_v^{\star}\circ S_ul=pS_v^{\star}s=-pq\,l$.  As a consequence, $S$
induces a non-degenerate metric if and only if the focal surfaces of
$f$ are non-degenerate, that is, are immersions into ${\cal Q}$.

Finally, we address the question: when is a given map
$S:M^2\to{\cal G}$ the conformal Gauss map of
some Legendre surface $f:M^2\to Z$?

Clearly, we have some necessary conditions:
first, the metric induced on $M$ by $S$ must have signature $(i,j)$ away
{}From critical points of $S$.
Secondly, with $(u,v)$ null coordinates for that metric, we have
$$\matrix{ S_u^{\star}\circ S_u=0 &{\rm and}& S_v\circ S_v^{\star}=0
\cr} \Eqno\Label\Necessary$$ (modulo interchanging the roles of $u$
and $v$).  Consider therefore $S:M\to{\cal G}$ satisfying these
conditions, and additionally assume that the induced metric is
non-degenerate.  In particular $S$ is an immersion.  From \Necessary,
we deduce ${\rm im}\,S_u\subset{\rm ker}\,S_u^{\star} =({\rm
im}\,S_u)^{\perp}$ so that ${\rm im}\,S_u$ is light-like and
1-dimensional (since $S_u\neq0$); similarly, ${\rm im}\,S_v^{\star}$
is light-like and 1-dimensional.  Thus, we obtain a candidate
$f=l\wedge s$, where $l,s:M\to{\cal L}^{\C}$ satisfy ${\rm
span}\{l\}={\rm im}\,S_v^{\star}$ and ${\rm span}\{s\}={\rm
im}\,S_u$.  Note that $f={\rm im}S_u\wedge{\rm im}S_v^{\star}$ and so
is real\Footnote{If we are in the situation of complex conjugate
parameters $(u,v)$, that is, $\bar S=S^{\perp}$ then we have
$S_v^{\star}=-\overline{S_u}$ so that ${\rm span}\{s\}={\rm
span}\{\bar l\}$.} since $dS$ is.

Since $l$ takes values in ${\rm im}\,S_v^{\star}\subset{\rm
ker}\,S_v=l^{\perp}\cap S$, we have $S_vl=0$, that is,
$l_v\in\Gamma(S)$; consequently, $l_v\in l^{\perp}\cap S$ which
implies $S_vl_v=0$, that is, $l_{vv}\in\Gamma(S)$.  We have therefore
established that $l,l_v,l_{vv}\in S$ and, similarly, $s,s_u,s_{uu}\in
S^{\perp}$.  We would like to show that $f$ is a Legendre immersion
which induces the same conformal structure on $M$ as $S$.  Sadly,
there are counterexamples to this assertion: one can construct $S$
with {\it constant\/} $f$.  However, as soon as $f$ is an immersion
the Legendre and conformality conditions amount to $l_u,s_v\in
l\wedge s=f$ and this is always true.  Indeed,
$S_ul=\pi^{\perp}l_u\in{\rm im}\,S_u={\rm span}\{s\}$ which gives us
half of the first assertion and it remains to prove $\pi l_u\in{\rm
span}\{l\}$.  At this point, we wheel out our non-degeneracy
assumption: since ${\rm im}\,S_u ={\rm span}\{s\}$ we get $S_u\circ
S_v^{\star}(s)=-\langle S_u,S_v\rangle s\neq0$ to find that
$S_v^{\star}s$ is a nonzero multiple of $l$.  Hence, $\pi l_u\in {\rm
span}\{l\}$ if and only if $\pi(S_v^{\star}s)_u\in{\rm span}\{l\}$.
The flat differentiation $d$ decomposes according to the bundle
decomposition $\C^6=S\oplus S^{\perp}$, $$\matrix{ d\sigma \hfill&=&
\nabla\sigma &+& dS\cdot\sigma &&{\rm for}& \sigma\in\Gamma(S)\ \cr
d\sigma^{\perp} &=& -dS^{\star}\cdot\sigma^{\perp} &+&
\nabla^{\perp}\sigma^{\perp} &&{\rm for}&
\sigma^{\perp}\in\Gamma(S^{\perp}), \cr }$$ and we obtain a Codazzi
equation $\nabla_vS_u^{\star}+S_v^{\star}\nabla^{\perp}_u
=\nabla_uS_v^{\star}+S_u^{\star}\nabla^{\perp}_v$ from the
$S^{\perp}\to S$ part of the (vanishing) curvature of $d$.  Apply
this to our distinguished section $s\in\Gamma(S^{\perp})$: since
$s,\nabla^{\perp}_vs\in s^{\perp}\cap S^{\perp}=({\rm
im}\,S_u)^{\perp} ={\rm ker}\,S_u^{\star}$, we have $\nabla_v
S^{\star}_u s-S^{\star}_u\nabla^{\perp}_vs=0$ so that the Codazzi
equation yields
$\nabla_u(S_v^{\star}s)=S_v^{\star}\nabla^{\perp}_us\in{\rm
span}\{l\}$.  This proves that $l_u\in l\wedge s$; similarly, we find
$s_v\in l\wedge s$.  To summarize:

\Thm[Theorem (Blaschke \Blaschke, \S93)]
{Let $S:M^2\to{\cal G}_{i,j}^{m,n}$ be an immersion which induces on $M$ a
metric of signature $(i,j)$ for which $(u,v)$ are null coordinates.
Then $S$ is the conformal Gauss map of a possibly degenerate Legendre map
$f:M^2\to Z$ if and only if\Footnote{It may happen that the same equations
additionally hold with the roles of $u$ and $v$ interchanged:
in this case, $f$ is the contact lift of a surface of Demoulin \Demoulin\
(or its Lie sphere geometry equivalent) where the congruence of Lie quadrics
only posesses two envelopes;
if only one of these additional conditions holds, we are in the case of
a Godeaux-Rozet surface \Godeaux, \Rozet\ (resp.\ its sphere geometry
equivalent).  See also \S122\ in \BolII, and \FerapontovD.}
$S_u^{\star}\circ S_u=0$ and $S_v\circ S_v^{\star}=0$.
In this case $f={\rm im}S_u\wedge{\rm im}S_v^{\star}$.}\Label\BlaschkeThm

\Section{The variational problem}
We now come to the main point of our considerations. Let $f:M^2\to Z$
be a Legendre map with non-degenerate induced conformal structure
and conformal Gauss map $S:M^2\to{\cal G}$.

\Def[Definition]{We define the Willmore energy of $f$ to be the harmonic
map energy of its conformal Gauss map $S$.  Thus, $$\matrix{
   W(f):=E(S)={1\over2}\int_M\langle dS,dS\rangle d{\rm vol}
        ={1\over2}\int_M\langle dS\wedge\star dS\rangle}$$
where $\star$ is the Hodge $\star$-operator on $M$ provided by the conformal
structure induced by $f$.

We say that $f$ is $W$-minimal if it extremizes $W$ with respect to
variations through Legendre maps.}

We shall see in Section~5 that the contact lift of an immersion $\f$
is $W$-minimal if and only if $\f$ is Lie or projectively minimal.

It is clear from the definition that $f$ extremizes $W$ as soon as
$S$ is harmonic.  Our mission is to prove the converse:

\Thm[Proposition]{If ${\partial\over\partial t}|_{t=0}W(f_t)=0$ for every
variation $f_t$ of $f_0=f$ through Legendre maps, then the conformal
Gauss map $S$ of $f$ is harmonic.}

This will require a little preparation.
\Comment
This requires some preparations: let $f_t$ be such a variation,
$S_t$ and $\star_t$ the corresponding conformal Gauss maps and Hodge
operators; further, let $D$ denote the pull back of the Levi-Civita
connection on ${\cal G}$ along $n_t$.
Then, ${\partial\over\partial t}W(f_t)
={1\over2}\int_M{\partial\over\partial t}\langle dS_t\wedge\star_tdS_t\rangle$;
for $t$-independent vector fields $x,y\in\Gamma(TM)$, we obtain
$$\matrix{ {\partial\over\partial t}\langle dS_t\wedge\star_tdS_t\rangle
   &=&\hfill \langle D_{\partial\over\partial t}dS_t(x),\star_tdS_t(y)\rangle
   +\langle dS_t(x),D_{\partial\over\partial t}(\star_tdS_t(y))\rangle\ \cr
   &&-\langle D_{\partial\over\partial t}dS_t(y),\star_tdS_t(x)\rangle
   -\langle dS_t(y),D_{\partial\over\partial t}(\star_tdS_t(x))\rangle. \cr }$$
As $D$ is torsion free, it follows that $$\matrix{
   D_{\partial\over\partial t}dS_t(x)
   =D_x{\partial\over\partial t}S_t+dS_t([{\partial\over\partial t},x])
   =D_x({\partial\over\partial t}S_t);}$$
and, writing $J_t$ for the adjoint of $\star_t$, that is,
$\star_tdf_t=df_t\circ J_t$, we find $$\matrix{
   D_{\partial\over\partial t}(\star_tdS_t(x))
   =D_{J_tx}({\partial\over\partial t}S_t)
      +dS_t([{\partial\over\partial t},J_tx])
   =D_{J_tx}({\partial\over\partial t}S_t)
      +dS_t(({\partial\over\partial t}J_t)x).}$$
Consequently, denoting $\dot{(..)}={\partial\over\partial t}|_{t=0}(..)$,
$$\matrix{
   {\partial\over\partial t}|_{t=0}\langle dS_t\wedge\star_tdS_t\rangle
   &=&\ \,\langle D\dot S\wedge\star dS\rangle
      +\langle dS\wedge\star D\dot S\rangle
      +\langle dS\wedge dS\circ\dot J\rangle \cr
   &=& 2\langle D\dot S\wedge\star dS\rangle
      +\langle dS\wedge dS\circ\dot J\rangle. \hfill\cr }$$
Next, we want to show that the second term vanishes,
$\langle dS\wedge dS\circ\dot J\rangle=0$.
Since $J_t^2\equiv\pm1$, we have $\dot JJ+J\dot J=0$.  Hence, $\dot J$
intertwines the eigenspaces of $J$: if $X_{\pm}$ denote eigenvectors of $J$,
$JX_{\pm}=\pm\varepsilon X_{\pm}$ with $\varepsilon=1$ in case $J^2=1$ and
$\varepsilon=i$ in case $J^2=-1$, then $\dot JX_{\pm}\parallel X_{\mp}$.
Moreover, the eigendirections of $J$ are isotropic for the conformal structure
showing that $\langle dS(X_{\pm}),dS(\dot JX_{\mp})\rangle=0$ since $S$ is
conformal, that is, $X_{\pm}$ are isotropic for the metric induced by $S$.

Finally, we compute $d\langle\dot S\wedge\star dS\rangle
=\langle D\dot S\wedge\star dS\rangle+\langle\dot S,d^D\star dS\rangle$;
so, integration by parts (Stoke's theorem) yields $$\matrix{
   {\partial\over\partial t}|_{t=0}W(f_t)
      &=&-\int_M\langle\dot S,d^D\!\star\!dS\rangle.\cr }$$
\EndComment

First, a standard computation in harmonic map theory (see for example
\BFLPP) gives $$\matrix{ {\partial\over\partial t}|_{t=0}W(f_t)
&=&-\int_M\langle\dot S,d^D\!\star\!dS\rangle +{1\over2}\int_M\langle
dS\wedge dS\circ\dot J\rangle,\cr }$$ where $D$ is the pull-back of
the Levi-Civita connection on ${\cal G}$ and $J_t$ is the adjoint of
$\star_t$, that is, $\star_tdf_t=df_t\circ J_t$.  We now show that
the second integrand vanishes: $\langle dS\wedge dS\circ\dot
J\rangle=0$.  Since $J_t^2\equiv\pm1$, we have $\dot JJ+J\dot J=0$.
Hence, $\dot J$ intertwines the eigenspaces of $J$: if $X_{\pm}$
denote eigenvectors of $J$, $JX_{\pm}=\pm\varepsilon X_{\pm}$ with
$\varepsilon=1$ in case $J^2=1$ and $\varepsilon=i$ in case $J^2=-1$,
then $\dot JX_{\pm}\parallel X_{\mp}$.  Moreover, the eigendirections
of $J$ are isotropic for the conformal structure showing that
$\langle dS(X_{\pm}),dS(\dot JX_{\mp})\rangle=0$ since $S$ is
conformal, that is, $X_{\pm}$ are isotropic for the metric induced by
$S$.

{\it Thus the conformality of $S$ allows us to ignore the variation in
conformal structure on $M$.}

In the case at hand, under the identification $S^{\ast}T{\cal G}
\cong{\rm Hom}(S,S^{\perp})$, the connection $D$ is induced by the
connections $\nabla$, $\nabla^{\perp}$ on $S$, $S^{\perp}$:
$D\tau=\nabla^{\perp}\circ\tau-\tau\circ\nabla.$
Fixing our null coordinates $(u,v)$ so that
$J{\partial\over\partial u}=\varepsilon{\partial\over\partial u}$ and
$J{\partial\over\partial v}=-\varepsilon{\partial\over\partial v}$
we compute $$\matrix{
   d^D\!\star\!dS({\partial\over\partial u},{\partial\over\partial v})
   &=\hfill& \varepsilon[-\nabla^{\perp}_u\circ S_v+S_v\circ\nabla_u
      -\nabla^{\perp}_v\circ S_u+S_u\circ\nabla_v] \cr
   &=\hfill&-2\varepsilon[\nabla^{\perp}_u\circ S_v-S_v\circ\nabla_u] \hfill\cr
   &=:& -2\varepsilon \tau_S \hfill\cr }$$
where the second equality follows from the Codazzi equation.
Thus, $$\matrix{{\partial\over\partial t}|_{t=0}W(f_t)
   &=& 2\varepsilon\int_M\langle\dot S,\tau_S\rangle\,du\wedge dv. \cr}$$

The key point now is that the mean curvature vector $\tau_S:S\to S^{\perp}$
of a conformal Gauss map has very restricted image --- it takes values in
(the pullback of) $L^{\star}$, the dual of the contact line bundle of $Z$:

\Def[Lemma]{${\rm im}\,\tau_S\subset{\rm span}\{s\}$, and
   ${\rm ker}\,\tau_S\supset l^{\perp}\cap S$.}

\Proof
{}From $\tau_S=\nabla^{\perp}_v\circ S_u-S_u\circ\nabla_v$ we learn
that ${\rm im}\,\tau_S\subset{\rm span}\{s\}$ since ${\rm im}\,S_u\subset
{\rm span}\{s\}$ and $\nabla^{\perp}_vs\in{\rm span}\{s\}$.  On the other hand,
the Codazzi equations yield $\tau_S=\nabla^{\perp}_u\circ S_v-S_v\circ\nabla_u$
so that the second claim follows since $l^{\perp}\cap S\subset{\rm ker}\,S_v$
and $l^{\perp}\cap S$ is $\nabla_u$-stable, since ${\rm span}\{l\}$ is.
\EndProof

Thus, $\tau_S\in{\rm Hom}(S/S\cap l^{\perp},{\rm span}\{s\})
\cong{\rm span}\{l\}\otimes{\rm span}\{s\}$ where we use the metric
to identify ${\rm span}\{l\}$ with $(S/S\cap l^{\perp})^{\star}$.

Using this, we have
$$(\tau_S^{\star}\circ\dot S)l={\rm tr}(\tau_S^{\star}\circ\dot S)\cdot
l =-\langle\dot S,\tau_S\rangle\cdot l$$ since $\tau_S^{\star}\in{\rm
Hom} (S^{\perp}/(s^{\perp}\cap S^{\perp}),{\rm span}\{l\})$, by the
lemma.  Moreover, fixing $\sigma\in\Gamma(S^{\perp})$ with $\langle
s,\sigma\rangle \equiv1$, we have $$\matrix{ \dot Sl=\langle s,\dot
Sl\rangle\sigma &{\rm mod}& s^{\perp}\cap S^{\perp}\cr }$$ whence
$(\tau_S^{\star}\circ\dot S)l=\langle s,\dot Sl\rangle
\tau_S^{\star}\sigma$.  Finally $\langle s,\dot Sl\rangle=\langle
s,\pi^{\perp}\dot l\rangle =\langle s,\dot l\rangle=\langle s,\dot
fl\rangle$ since $\dot fl=\dot l$ mod $f$.  Thus $$\matrix{ \langle
s,\dot Sl\rangle = -\vartheta(\dot f)(s,l) \cr}$$ where $\vartheta$
is (the pullback of) our contact form so that
$(\tau_S^{\star}\circ\dot S)l=-\vartheta(\dot
f)(s,l)\tau_S^{\star}\sigma$.  Otherwise said
$(\tau_S^{\star}\circ\dot S)l=-\vartheta(\dot f)
(s,\tau_S^{\star}\sigma)l$ and taking a trace gives $$\matrix{
\langle\dot S,\tau_S\rangle = \vartheta(\dot
f)(s,\tau_S^{\star}\sigma).\cr }$$

Now let $\gamma$ be an arbitrary compactly supported section of $L$.  It is
well known\Footnote{Indeed, take any vector field $X$ with $\vartheta(X)
=\gamma$ and observe that $\Gamma({\cal D})\ni Y\mapsto\vartheta([X,Y])\in L$
is tensorial.  Now use the nondegeneracy of $\vartheta([.,.])$ on ${\cal D}$
to find $Y\in\Gamma({\cal D})$ with $\vartheta([X,.])=\vartheta([Y,.])$ on
$\Gamma({\cal D})$ and set $X_{\gamma}=X-Y$.} \Liebermann\ that $\gamma$
generates an infinitesimal contactomorphism $X_{\gamma}$ on $Z$ with
$\vartheta(X_{\gamma})=\gamma$.
Let $\Phi_t:Z\to Z$ be the flow by contactomorphisms of $X_{\gamma}$ and
set $f_t=\Phi_t\circ f$.  Then $f_t$ is a variation of $f$ through Legendre
maps with $\dot f=X_{\gamma}\circ f$ and we have $$\matrix{
   \langle\dot S,\tau_S\rangle = \gamma|_f(s,\tau_S^{\star}\sigma) .\cr}$$
Thus, if $\D{}{t}|_{t=0}W(f_t)=0$, we see that $\tau_S^{\star}\sigma$
vanishes whence $\tau_S^{\star}$ and so $\tau_S$ vanish by the
lemma.  We conclude that $S$ is harmonic.
To summarize:

\Thm[Theorem]{A Legendre surface $f:M\to Z$ extremizes $W$ with
respect to Legendre variations if and only if
its conformal Gauss map $S:M\to{\cal G}$ is harmonic.}

\Comment
Now, let $\gamma$ be {\it any\/} skew scalar product on the tautological
bundle over $Z$ (remember that those form a line bundle $L$ over $Z$, where
the contact form takes values).  We intend to show how to construct a variation
$f_t$ through Legendre surfaces with $\gamma|_f=\vartheta(\dot f)$ from such
$\gamma$, chosen arbitrarily enough: this will establish the our theorem since
it will force $\tau_S^{\star}$ to vanish, $\tau_S^{\star}\hat s=0$.

Thus, choose some vector field $X\in\Gamma(Z)$ with $\vartheta(X)=\gamma$.
For $Y\in\Gamma({\cal D})$ and $g\in C^{\infty}(Z)$ we compute
$\vartheta([X,gY])=g\vartheta([X,Y])+(Xg)\vartheta(Y)=g\vartheta([X,Y])$,
so that $\vartheta([X,.]):\Gamma({\cal D})\to L$ behaves tensorial.
However, $-d\vartheta=\vartheta([.,.])$ is a non-degenerate 2-form
on ${\cal D}$; so, there is a (unique) vector field $Y\in\Gamma({\cal D})$
such that $\vartheta([X,.])=\vartheta([Y,.])$.  Set $X_{\gamma}:=X-Y$
to conclude $\vartheta([X_{\gamma},.])\equiv0$.  In other words,
$[X_{\gamma},.]:\Gamma({\cal D})\to\Gamma({\cal D})$, that is, $X_{\gamma}$
is an infinitesimal symplectomorphism\Footnote{In fact, this gives an
isomorphism between $\Gamma(L)$ and the Lie algebra of infinitesimal
symplectomorphisms on $Z$.}.

The reason here is that, for any 1-parameter group $\Phi_t$ of
diffeomorphisms, the derivative of its corresponding vector field
$X=\dot\Phi$ is given by $\partial_t|_{t=0}(d\Phi_t(Y))=[X,Y]$
--- actually, this is the way the Lie derivative is often introduced
(cf.Gallot-Hulin-Lafontaine, or Sternberg).  And, if the derivative of
such a 1-parameter group preserves the contact distribution, then
$\Phi_t$ does...

Now, take $\gamma$ (and, therefore, $X_{\gamma}$) to be compactly
supported, and let $\Phi_t:Z\to Z$ be the associated flow of
contactomorphisms.
Set $f_t:=\Phi_t\circ f$.  This gives a variation of $f$ through
Legendre maps, with $\dot f=X_{\gamma}\circ f$.  Consequently, we have
$$\matrix{
   \langle\dot S,\tau_S\rangle = -\gamma|_f(\tau_S^{\star}\hat s,s). }$$
Since $\gamma\in\Gamma(L)$ was an arbitrarily compactly supported skew
product on the tautological bundle over $Z$, the theorem follows.
\EndComment

\Section{Implications}
We now explore the implications of our analysis and consider each of the
signatures $(m,n)$ in turn.  In this way, we obtain results on Lie and
projectively minimal surfaces as promised in the introduction.

\SubSection{Projectively minimal surfaces}
Consider first the case $(m,n)=(3,3)$.
Here we are dealing with projective differential geometry.
For this, view $\RP^3$ as an ${\rm SL}(4,\R)$-space: thus we give
$\R^4$ a fixed volume form ${\rm vol}\in\Lambda^4(\R^4)^{\star}$ and set
$\RP^3=\{\R x\,|\,x\in\R^4\setminus\{0\}\}$.
The 6-dimensional space $\Lambda^2(\R^4)$ gets a metric of signature $(3,3)$
by $\langle v,w\rangle={\rm vol}(v\wedge w)$ for which the action of
${\rm SL}(4,\R)$ is clearly isometric.
This gives a double covering ${\rm SL}(4,\R)\to{\rm O}(3,3)$.
Henceforth, we write $\Lambda^2(\R^4)=\R^{3,3}$.
Moreover, $l\in\R^{3,3}$ satisfies the Pl\"ucker relation
$\langle l,l\rangle=0$ if and only if $l$ is decomposible:
$l=x\wedge y$ for some $x,y\in\R^4$.
This yields a diffeomorphism, the {\it Klein correspondence\/},
between lines in $\RP^3$ and the Pl\"ucker quadric ${\cal Q}$:
${\rm span}\{x,y\}\mapsto{\rm span}\{x\wedge y\}$.
Further, if $l_1,l_2\in{\cal Q}$ satisfy $\langle l_1,l_2\rangle=0$ then
there is an $x\in\R^4$ such that $l_1=x\wedge y_1$ and $l_2=x\wedge y_2$.
That is, orthogonal points in ${\cal Q}$ correspond to intersecting lines
in $\RP^3$.
Moreover, ${\rm span}\{l_1,l_2\}=\{x\wedge(ay_1+by_2)\,|\,a,b\in\R\}$
determines a plane ${\rm span}\{x,y_1,y_2\}$ in $\RP^3$.
Hence, we can identify the space $Z$ of null lines in ${\cal Q}$ with the
manifold $\{(p,P)\in\RP^3\times(\RP^3)^{\ast}\,|\,p\in P\}$ of contact
elements to $\RP^3$.

Let ${\fk f}:M^2\to\RP^3$ be an immersion and contemplate its contact
lift $f:M^2\to Z$ defined by $f={\rm span}\{{\fk f}\wedge d{\fk
f}(X)\,|\,X\in TM\}$.

$f$ is an immersion: indeed, if not, there is some $X\in T_xM$ with
$d_Xl\in f(x)$ for any local section $l$ of $f$.
In particular, ${\fk f}(x)\wedge dl(X)=0$.
However, any such $l$ is of the form ${\fk f}\wedge d_Y{\fk f}$
for some local vector field $Y\in\Gamma(TM)$ so that
$0={\fk f}\wedge d_X({\fk f}\wedge d_Y{\fk f})
  ={\fk f}\wedge d_X{\fk f}\wedge d_Y{\fk f}$
for all $Y\in T_xM$.
Consequently,
$ d_X{\fk f}\in\bigcap_{Y}{\rm span}\{{\fk f}(x), d_Y{\fk f}\}
={\rm span}\{{\fk f}(x)\}$, a contradiction to the assumption that ${\fk f}$
is an immersion.

$f$ is Legendre: fixing a local basis $(X,Y)$ in $TM$, set
$l:={\fk f}\wedge d_X{\fk f}$ and $s:={\fk f}\wedge d_Y{\fk f}$,
so that $f={\rm span}\{l,s\}$.
Then, $\langle dl,s\rangle{\rm vol}
  =d{\fk f}\wedge d_X{\fk f}\wedge{\fk f}\wedge d_Y{\fk f}=0$
since $d{\fk f}$ certainly takes values in
${\rm span}\{{\fk f}, d_X{\fk f}, d_Y{\fk f}\}$.

Now, we identify the conformal structure induced by $f$: first note that
the conformal class of the second fundamental form of ${\fk f}$ is a
projective invariant.  Indeed, let ${\fk f}^{\star}:M^2\to(\RP^3)^{\star}$
denote the dual surface of ${\fk f}$ and contemplate $TM\ni X,Y\mapsto
{\fk f}^{\star}( d_X d_Y{\fk f})$: this is tensorial and scales
with ${\fk f}$ and ${\fk f}^{\star}$.
We fix a basis $(X,Y)$ in $TM$ and choose ${\fk n}$ such that
${\rm vol}({\fk f}, d_X{\fk f}, d_Y{\fk f},{\fk n})\equiv1$;
we adjust the scaling of ${\fk f}^{\star}$ so that
${\fk f}^{\star}({\fk n})\equiv1$.
Then,
${\fk f}\wedge d_X{\fk f}$ and ${\fk f}\wedge d_Y{\fk f}$
form a basis of $f$ while
${\fk f}\wedge{\fk n}$ and $ d_X{\fk f}\wedge d_Y{\fk f}$ mod $f$
form a basis of $f^{\perp}/f$.
Now, for $Z=aX+bY\in TM$ we find $$\matrix{
   df(Z)\cdot{\fk f}\wedge d_X{\fk f}
      &=&\hfill  d_Z{\fk f}\wedge d_X{\fk f}
      &+& {\fk f}\wedge d_Z d_X{\fk f} &{\rm mod}& f\cr
   &=& -b\, d_X{\fk f}\wedge d_Y{\fk f} 
      &+& {\fk f}^{\star}( d_Z d_X{\fk f})
         {\fk f}\wedge{\fk n} &{\rm mod}& f\cr
   df(Z)\cdot{\fk f}\wedge d_Y{\fk f}
      &=&\hfill a\, d_X{\fk f}\wedge d_Y{\fk f}
      &+& {\fk f}^{\star}( d_Z d_Y{\fk f})
         {\fk f}\wedge{\fk n} &{\rm mod}& f\cr}$$
so that the determinant of $df(Z)$
--- that gives, by definition, the conformal structure ---
is $-a{\fk f}^{\star}( d_Z d_X{\fk f})-b{\fk f}^{\star}
( d_Z d_X{\fk f})=-{\fk f}^{\star}( d_Z^2{\fk f})$.
Consequently, the conformal class on $M$ induced by $f$ is that of the
second fundamental form of ${\fk f}$.
It follows that the null directions of the conformal structure are precisely
the asymptotic directions of ${\fk f}$.

Finally, let $u,v$ be (possibly complex conjugate) asymptotic coordinates
for ${\fk f}$, and define $l:={\fk f}\wedge{\fk f}_u$ and
$s:={\fk f}\wedge{\fk f}_v$.
Then, $f={\rm span}\{l,s\}$, and $l$, $s$ are by construction the line
congruences tangent to the asymptotic directions of ${\fk f}$.
We have ${\fk f}_{uu},{\fk f}_{vv}\in{\rm span}\{{\fk f},{\fk f}_u,{\fk f}_v\}$
so that $$\matrix{
   {\fk f}_{uu} &=& \ast\,{\fk f}_u &+& p\,{\fk f}_v &+& \ast\,{\fk f}, \cr
   {\fk f}_{vv} &=& q\,{\fk f}_u &+& \ast\,{\fk f}_v &+& \ast\,{\fk f}. \cr
   }$$
with suitable functions $p,q$.  Now, $$\matrix{
   l_u &=& ({\fk f}\wedge{\fk f}_u)_u
       &=& {\fk f}_u\wedge{\fk f}_u + {\fk f}\wedge{\fk f}_{uu}
       &=& \ast\,l &+& ps, \cr
   s_v &=& ({\fk f}\wedge{\fk f}_v)_v
       &=& {\fk f}_v\wedge{\fk f}_v + {\fk f}\wedge{\fk f}_{vv}
       &=& ql &+& \ast\,s. \cr
   }$$
{}From this we see that $l$ and $s$ are the focal surfaces of the line
congruence $f$ in ${\cal Q}$, with $u$ and $v$ the corresponding conjugate
parameters, and so we are in the situation of our main analysis.
Moreover, from \WillmoreDensity, we learn that the Willmore energy of $f$
is given by $$\matrix{
   W(f)=\int_M\langle S_u,S_v\rangle du\wedge dv=\int_Mpq\,du\wedge dv.}$$
This coincides with the projectively minimal Lagrangian $L_{\rm proj}({\fk f})$
for the immersion ${\fk f}$ described in the introduction.
All our constructions are palpably ${\rm SL}(4,\R)$-invariant and we conclude,
with \ThomsenA, that $L_{\rm proj}({\fk f})$ is a projectively invariant
functional.
Moreover, as ${\fk f}$ varies through immersions, its contact lift varies
through Legendre immersions and conversely so that ${\fk f}$ is projectively
minimal if and only if its contact lift $f$ is $W$-minimal.

What is the geometry of the conformal Gauss map $S$ of $f$?
We have already seen that any $S(x)$, $x\in M$, gives rise to a Legendre
submanifold $Z_{S(x)}$ having first order contact with $f$ at $f(x)$.
We claim that $Z_{S(x)}$ is the contact lift of a quadric
$Q_{S(x)}\subset\RP^3$ which therefore has second order contact
with the underlying immersion ${\fk f}:M\to\RP^3$ at ${\fk f}(x)$.

A quadric $Q$ in $\RP^3$ is the null cone of an inner product, also
called $Q$, on $\R^4$ which is unique up to homothety.  The quadric has
a conformal structure (given by the second fundamental form) of signature
$(i,j)$ which is $(1,1)$ when the metric $Q$ has signature $(2,2)$ and is
$(2,0)$ when the metric has Lorentz signature.  Without loss of generality,
we can take ${\rm vol}_Q={\rm vol}$ so that the Hodge operator
$\star_Q:\Lambda^2\R^4\to\Lambda^2\R^4$ is given by $$\matrix{
   {\rm vol}(v\wedge\star_Qw) = Q(v,w). \cr}$$
Now $\star_Q$ is clearly a symmetric endomorphism of $\Lambda^2\R^4=\R^{3,3}$
(remember that the $(3,3)$ metric on $\Lambda^2\R^4$ is provided by
${\rm vol}$) and $\star_Q^2=\pm1$ according to whether the metric $Q$ has
signature $(2,2)$ or Lorentz signature.  Otherwise said, we have defined
a map $Q\mapsto\star_Q$ from the space of quadrics with $(i,j)$ conformal
structure to ${\cal G}^{3,3}_{i,j}$.
Both domain and co-domain of this map are homogeneous ${\rm SL}(4,\R)$-spaces
while this map is plainly ${\rm SL}(4,\R)$-equivariant and so is a surjection.
In fact, it injects also so that a quadric is determined by its $\star$-operator
on 2-vectors.  Indeed, as is well known, we have

\Def[Lemma]
{\it Let $l=x\wedge y\in\Lambda^2(\R^4)$ be decomposable.
Then $l$ is an eigenvector of $\star_Q$, $\star_Ql=\pm\varepsilon l$, if and
only if $l$ is a line on the quadric given by $Q$, that is, $Q$ vanishes on
the 2-plane in $\R^4$ spanned by $x$ and $y$.}

With this in hand,  we see that the two families of generators of a quadric
$Q$ are given by the null (for the $(3,3)$ metric) eigenvectors of $\star_Q$.
{}From this, it is clear that $\star_Q$ determines the quadric $Q$ and, moreover,
for $S\in{\cal G}^{3,3}_{i,j}$, $$\matrix{
   Z_S &=& \{{\rm span}\{l_+,l_-\}\,|\,\star_Sl_{\pm}=\pm\varepsilon l_{\pm}\}
   \cr}$$
comprises the span of pairs of generators, one from each family or,
equivalently, the contact element given by the intersection of these
lines, together with their span in $T\RP^3$.  Otherwise said, $Z_S$ is
the contact lift of the corresponding quadric $Q_S$.
To summarize, the conformal Gauss map $S$ of $f$ is (the contact lift of)
a congruence of quadrics in $\RP^3$ having second order contact with
the surface ${\fk f}$.

Classically, these quadrics are known as the ``Lie quadrics'' (see \Blaschke,
\BolII, \ThomsenA, compare \LieBrief) so that the conformal Gauss map of our
Legendre surface $f$ is the congruence of Lie quadrics of ${\fk f}$.
Putting this all together, we have:

\Thm[Theorem]
{\it A surface ${\fk f}:M^2\to\RP^3$ is projectively minimal if and only if
its congruence of Lie quadrics\Footnote{Here, $(i,j)$ is the signature of
the second fundamental form of ${\fk f}$ and so is $(2,0)$ for convex
surfaces and $(1,1)$ for negatively curved ${\fk f}$ --- our non-degeneracy
assumption on $f$ is precisely the non-degeneracy of this second fundamental
form.} $S:M^2\to{\cal G}^{3,3}_{i,j}$ is a harmonic map.}

We complete this circle of ideas by briefly discussing how to recover
an immersion ${\fk f}$ from a Legendre surface $f$.
Clearly, {\it any\/} surface $f={\rm span}\{l,s\}$ gives a map
${\fk f}:M^2\to\RP^3$ as $\langle l,s\rangle=0$ implies
$l={\fk f}\wedge{\fk g}_1$ and $s={\fk f}\wedge{\fk g}_2$
and then the contact condition gives us
$0=d({\fk f}\wedge{\fk g}_1)\wedge({\fk f}\wedge{\fk g}_1)
   =d{\fk f}\wedge{\fk g}_1\wedge{\fk f}\wedge{\fk g}_2$,
that is, $d{\fk f}$ takes values in
${\rm span}\{{\fk f},{\fk g}_1,{\fk g}_2\}$.
The main issue is whether ${\fk f}$ is an immersion: this is
the case as soon as the conformal structure induced by $f$
is non-degenerate.
To see this, introduce ${\fk n}$ and ${\fk f}^{\star}$ as above
so that ${\rm vol}({\fk f},{\fk g}_1,{\fk g}_2,{\fk n})\equiv1$,
and ${\fk f}^{\star}$ annihilates ${\fk f}$, ${\fk g}_1$, ${\fk g}_2$,
and has ${\fk f}^{\star}({\fk n})\equiv1$.
Now suppose ${\fk f}$ is not immersed at some point $x\in M$.
Then, there are two cases to consider:
if $d{\fk f}\equiv0$ mod ${\fk f}$ at $x$, then
$df(Z)\cdot{\fk f}\wedge{\fk g}_i
   ={\fk f}^{\star}( d_Z{\fk g}_i){\fk f}\wedge{\fk n}$ mod $f$
for any $Z\in T_xM$, and the induced conformal structure degenerates
completely.
If, on the other hand, $d{\fk f}\not\equiv0$ mod ${\fk f}$ at $x$, we can
assume without loss of generality that there is a basis $(X,Y)$ in $T_xM$
such that $ d_X{\fk f}\parallel{\fk f}$ and $ d_Y{\fk f}={\fk g}_2$
--- implying that $ d_X({\fk f}\wedge{\fk g}_2)\equiv0$ mod $f$.
For $Z=aX+bY$ we then find $$\matrix{
   df(Z)\cdot{\fk f}\wedge{\fk g}_1
      &=&\hfill -b{\fk g}_1\wedge{\fk g}_2
      &+&\hfill {\fk f}^{\star}( d_Z{\fk g}_1){\fk f}\wedge{\fk n}
      &{\rm mod}& f, \cr
   df(Z)\cdot{\fk f}\wedge{\fk g}_2
      &=&
      &&\hfill b{\fk f}^{\star}( d_Y{\fk g}_2){\fk f}\wedge{\fk n}
      &{\rm mod}& f. \cr }$$
Hence, the induced conformal structure is given by
   $aX+bY=Z\mapsto-b^2{\fk f}^{\star}( d_Y{\fk g}_2)$
giving $X\perp T_xM$ so that, in this case also, the conformal structure
degenerates.
As a consequence, any Legendre surface $f:M^2\to Z$ with a non-degenerate
conformal structure gives rise to an immersion ${\fk f}:M^2\to\RP^3$ into
projective $3$-space.

\SubSection{Lie minimal surfaces}
Here, we skip most of the setup --- all the background material can
be found in the modern and readable introduction by Cecil \Cecil.
The key points, however, are that the points of the quadric
${\cal Q}\subset\P\R^{4,2}$ parametrize oriented spheres in
$S^3\cong\R^3\cup\{\infty\}$ and that two such points are orthogonal
if and only if the corresponding spheres are in oriented contact.
Thus we may identify $Z$ with a space of contact elements, this time
the contact elements of $S^3$.

Let ${\fk f}:M^2\to\R^3$ be an immersion with Gauss map ${\fk
n}:M^2\to S^2$.  Let $(v_{-1},v_0,v_1,v_2,v_3,v_{\infty})$ denote a
basis of $\R^{4+2}$ such that $(v_1,v_2,v_3)$ form an orthonormal
basis for an $\R^{3}$, $v_{-1}$ is an orthogonal time-like basis
vector, and $(v_0,v_{\infty})$ are isotropic, spanning the remaining
orthogonal $\R^{1+1}$, with scalar product $\langle
v_0,v_{\infty}\rangle=-{1\over2}$.  Then, define $\varphi:=v_0+{\fk
f}+{\fk f}^2v_{\infty}$ to be the stereographic projection of ${\fk
f}$ in $S^3$, $\nu:=v_{-1}+{\fk n}+2{\fk n}\cdot{\fk f}\,v_{\infty}$
its tangent plane map (here, ``$\cdot$'' denotes the scalar product
in $\R^3$) and let $f=\varphi\wedge\nu$ be the corresponding contact
lift.  Choose (local) curvature line coordinates $(u,v)$ so that
$0={\fk n}_u+\kappa_1{\fk f}_u$ and $0={\fk n}_v+\kappa_2{\fk f}_v$
where $\kappa_i$ are the principal curvatures of ${\fk f}$.  Further,
define the curvature spheres $l,s:M^2\to{\cal Q}$ by
$l=\nu+\kappa_1\varphi$ and $s=\nu+\kappa_2\varphi$.  Clearly,
$f=l\wedge s$, while $$\matrix{ l_u=(\partial_u\kappa_1)\varphi\in f
&{\rm and}& s_v=(\partial_v\kappa_2)\varphi\in f. \cr
}\Eqno\Label\CurvatureSpheres$$ Thus, the curvature spheres are the
focal surfaces of $f$ and the curvature line coordinates $(u,v)$ the
corresponding conjugate coordinates so that, once again, our analysis
applies.  Finally, rearranging \CurvatureSpheres, we have $$\matrix{
l_u &=& { \partial_u\kappa_1\over\kappa_1-\kappa_2}(l-s) &{\rm and}&
s_v &=& { \partial_v\kappa_2\over\kappa_1-\kappa_2}(l-s), \cr }$$ and
then \WillmoreDensity\ gives $$\matrix{ \langle S_u,S_v\rangle\,
du\wedge dv &=& -{ \partial_u\kappa_1
\partial_v\kappa_2\over(\kappa_1-\kappa_2)^2} du\wedge dv. \cr}$$
This coincides with the Lie minimal Lagrangian $L_{\rm Lie}({\fk f})$
for the immersion ${\fk f}$ as described in the introduction.  As all
our constructions were Lie-invariant we conclude that the functional
$L_{\rm Lie}({\fk f})$ is invariant under Lie sphere transformations.
Moreover, as ${\fk f}$ varies through immersions, its contact lift
varies through Legendre immersions and conversely\Footnote{Note that
the condition on ${\fk f}$ to be an immersion is an open condition.}
so that ${\fk f}$ is Lie minimal if and only if its contact lift $f$
is $W$-minimal.

Just as in the projective case, the conformal Gauss map $S$ of a Legendre
surface $f$ provides a congruence of ``simple'' surfaces having second
order contact with the underlying surface:
in the case at hand, the conformal Gauss map defines a congruence of
Dupin cyclides.  This becomes clear when realizing that the spheres in
each $S(x)$ and $S^{\perp}(x)$ are the principal spheres of $Z_{S(x)}$.
Namely, the principal spheres of $Z_{S(x)}$ are constant along the
corresponding curvature lines so that the surface has two families
of circular curvature lines (is a channel surface in two ways).
As in the projective case, the Dupin cyclides of the conformal Gauss
map are the Lie cyclides \LieBrief\ of the Legendre surface $f$.  Thus:

\Thm[Theorem]
{\it A surface ${\fk f}:M^2\to\R^3\subset S^3$ is Lie minimal if and only if
its congruence of Lie cyclides $S:M^2\to{\cal G}^{4,2}_{1,1}$ is a harmonic
map.}

We want to conclude this section with a remark on the recovery of a surface
in $\R^3$ or $S^3$ from a Legendre surface $f:M^2\to Z^{4,2}$:
in contrast to the projective picture, not every Legendre surface gives
rise to an immersion ${\fk f}:M^2\to\R^3$.  An easy example is given by
the horn and spindle cyclides\Footnote{These are M\"obius transformations
of circular cylinders and circular cones, respectively.} where the
Legendre surface is immersed but an induced map ${\fk f}:M^2\to\R^3$
becomes singular at certain points.  However, as all Dupin cyclides
are Lie equivalent, one can find an immersed map ${\fk f}:M^2\to\R^3$
with contact lift $f$, parametrizing a ring cyclide (torus) in this
example.
For a more comprehensive discussion see \Cecil.

\Section{Applications}
We conclude our discussion by recalling, in abbreviated form, the main points
of the integrable systems theory of harmonic maps and indicating how these
apply to $W$-minimal surfaces.

\SubSection{Harmonic maps and moving frames}
Let $N$ be a pseudo-Riemannian symmetric $G$-space for some Lie group $G$.
Fix a base-point $o\in N$ with stabilizer $K$ so that $N\cong G/K$ and the
involution at $o$ induces a symmetric decomposition ${\fk g}={\fk k}\oplus
{\fk p}$.
The coset projection $p:G\to N$, $p(g)=g\cdot o$, is a principal $K$-bundle.

Now let $M$ be a surface equipped with a conformal structure of signature
$(i,j)$ and contemplate maps $\varphi:M\to N$.
A {\it frame\/} of $\varphi$ is a map $F:M\to G$ such that $p\circ F=\varphi$.
Let $F$ be such a frame and consider its Maurer-Cartan form $\alpha=F^{-1}dF$:
a 1-form on $M$ with values in ${\fk g}$.
Write $\alpha=\alpha_{\fk k}+\alpha_{\fk p}$ according to the symmetric
decomposition and further write $$\matrix{
   \alpha_{\fk p} = \alpha_{\fk p}' + \alpha_{\fk p}'' \cr}$$
where $\alpha_{\fk p}'$ and $\alpha_{\fk p}''$ are the components along
the null directions of the conformal structure on $M$.
Thus, if $(u,v)$ are an oriented choice of null coordinates,
$\alpha_{\fk p}'=\alpha_{\fk p}(\D{}{u})du$.
Note that $\alpha_{\fk p}'$, $\alpha_{\fk p}''$ are real, respectively
complex conjugate according as $(i,j)$ is $(1,1)$ or $(2,0)$.

Now introduce a spectral parameter $\lambda\in \C^{\times}$ and
a family of ${\fk g}^{\C}$-valued 1-forms by $$\matrix{
   \alpha_{\lambda} := \alpha_{\fk k}
      + \lambda\alpha_{\fk p}' + \lambda^{-1}\alpha_{\fk p}'' \cr}
   \Eqno\Label\SpectralFamily$$
The key observation (see \Pohlmeyer, \ZSh, \Uhlenbeck) is that $\varphi$
is harmonic if and only if each $\alpha_{\lambda}$ is flat: $$\matrix{
   d\alpha_{\lambda}+{1\over2}[\alpha_{\lambda}\wedge\alpha_{\lambda}]=0.
   \cr}\Eqno\Label\MaurerCartan$$
Conversely, given $\alpha_{\lambda}$ of the form \SpectralFamily\
satisfying \MaurerCartan\ for all $\lambda\in\C^{\times}$, we
may (locally) integrate to find $F_{\lambda}:M\to G^{\C}$ with
$F_{\lambda}^{-1}dF_{\lambda}=\alpha_{\lambda}$.
In particular, taking $\lambda\in\R^{\times}$ when $(i,j)=(1,1)$
and $\lambda\in S^1$ when $(i,j)=(2,0)$, we see that $\alpha_{\lambda}$
is ${\fk g}$-valued so that we may take $F_{\lambda}:M\to G$ and so obtain
a 1-parameter family of maps $\varphi_{\lambda}=p\circ F_{\lambda}:M\to N$
with $\varphi_1=\varphi$.

In the case at hand, we take $N={\cal G}^{m,n}_{i,j}$ with base-point $S_o$.
Then ${\fk g}={\fk o}(m,n)$ and $$\matrix{
   {\fk k}^{\C} &=& \{X:\C^6\to\C^6\,|\,X^{\star}=-X, XS_o\subset S_o\},\hfill\cr
   {\fk p}^{\C} &=& \{X:\C^6\to\C^6\,|\,X^{\star}=-X, XS_o\subset S_o^{\perp}\}.
      \cr}$$
Given $S:M^2\to{\cal G}^{m,n}_{i,j}$ with frame $F$, we have $$\matrix{
   F\alpha_{\fk p}'(\D{}{u})F^{-1} = S_u - S_u^{\star} &{\rm and}&
   F\alpha_{\fk p}''(\D{}{v})F^{-1} = S_v - S_v^{\star}. \cr
   }\Eqno\Label\StructureS$$
With this in hand, we offer some sample applications.

\SubSection{Spectral deformation}
It is easy to see that $\varphi_{\lambda}:M\to N$ introduced above are all
harmonic \Uhlenbeck\ so that (locally) we obtain a 1-parameter family of
harmonic maps from a given one.
In our setting, if $S$ is the conformal Gauss map of a $W$-minimal Legendre
surface, then $S_u^{\star}\circ S_u=S_v\circ S_v^{\star}=0$.
{}From \StructureS, we see that this condition amounts to $$\matrix{
   \alpha_{\fk p}'(\D{}{u})\circ\alpha_{\fk p}'(\D{}{u})\,S_o = 0 &{\rm and}&
   \alpha_{\fk p}''(\D{}{v})\circ\alpha_{\fk p}''(\D{}{v})\,S_o^{\perp} = 0.\cr
   }\Eqno\Label\HarmonicS$$
The lift $F_{\lambda}$ of $S_{\lambda}$ has
$(\alpha_{\lambda})_{\fk p}'=\lambda\alpha_{\fk p}'$
$(\alpha_{\lambda})_{\fk p}''=\lambda^{-1}\alpha_{\fk p}''$
which clearly satisfies \HarmonicS\ also so that, by Theorem \BlaschkeThm,
we see that the spectral deformation $S\mapsto S_{\lambda}$ preserves the
property of being a conformal Gauss map.
We therefore conclude (as has Ferapontov \FerapontovP\ and \FerapontovL\ by
different methods) that

\Thm[Theorem]
{A $W$-minimal surface $f:M^2\to Z$ gives (locally) rise to a 1-parameter
family of $W$-minimal surfaces $f_{\lambda}:M^2\to Z$ with $\lambda\in\R
^{\times}$ or $S^1$ according to the signature of the conformal
structure induced by $f$.}

\SubSection{Dressing transformations and B\"acklund transforms}
The lifts $F_{\lambda}$, $\lambda\in\C^{\times}$, patch together to
give a map of $M$ into a loop group.
There is a well developed theory of dressing actions where a point-wise
action of a complementary loop group induces an action on such maps and so,
eventually, on the underlying harmonic maps (see, for example, \BP).
In general, such an action requires solution of a Riemann-Hilbert problem
but for certain special elements of the complementary loop group, the
{\it simple factors\/}, the action is explicitly computable and gives
rise to B\"acklund transformations.
For a careful account of these ideas see \TerngUhlenbeck.

All we want to say here is that the dressing action preserves the class
of conformal Gauss maps.  For this, the only fact we need is that if
$\hat\varphi$ arises from $\varphi$ by a dressing transformation then
there are frames $F$, $\hat F$ with $$\matrix{
   \hat\alpha_{\fk p}' = k_+^{}\alpha_{\fk p}'k_+^{-1} &{\rm and}&
   \hat\alpha_{\fk p}'' = k_-^{}\alpha_{\fk p}''k_-^{-1} \cr }$$
where $k_{\pm}:M\to K^{\C}$.  Clearly the requirement \HarmonicS\ is
invariant under conjugation by elements of $K^{\C}$ and we conclude:

\Thm[Theorem]
{The dressing action on harmonic maps $M\to{\cal G}^{m,n}_{i,j}$ preserves
the class of conformal Gauss maps and so there is an induced action of a
loop group on $W$-minimal surfaces.}

In particular, dressing by simple factors provides B\"acklund
transforms of Lie and projectively minimal surfaces.  We will return
to this topic elsewhere.

\SubSection{A duality between harmonic maps}
The celebrated duality between Riemannian symmetric spaces of compact
and non-compact type \Helgason\ extends to arbitrary symmetric spaces:
if ${\fk g}={\fk k}\oplus{\fk p}$ is a symmetric decomposition then
$\hat{\fk g}={\fk k}\oplus\sqrt{-1}\,{\fk p}$ is also a symmetric
decomposition of a second Lie algebra giving rise to a second symmetric
space $\hat N$.

Now suppose that $\varphi:M^2\to N$ is a harmonic map of a surface with
$(1,1)$ conformal structure.
Then the spectral deformation arises by integrating $\alpha_{\lambda}$
for $\lambda\in\R^{\times}$.
However, for $\lambda\in\sqrt{-1}\,\R^{\times}$, $\alpha_{\lambda}$
is clearly $\hat{\fk g}$-valued so that $F_{\lambda}$ may be taken\Footnote{%
By adjusting constants of integration.} to be $\hat G$-valued and then
$\varphi_{\lambda}=p\circ F_{\lambda}$ is a harmonic map $M\to\hat N$.
We therefore have:

\Thm[Proposition]
{Let $M$ be a surface with $(1,1)$ conformal structure.
Then there is a (local) bijective correspondence between harmonic maps
$\varphi:M\to N$ and $\hat\varphi:M\to\hat N$ modulo isometries.}

In the case at hand, we have a duality between ${\cal G}^{3,3}_{1,1}$ and
${\cal G}^{4,2}_{1,1}$: the duality mechanism is implemented by taking
$\R^{3,3}=S_o\oplus S_o^{\perp}$ and then setting
$\R^{4,2}=S_o\oplus\sqrt{-1}\,S_o^{\perp}$.
Just as in Section~6.2, the duality preserves conformal
Gauss maps and we arrive at a conceptual explanation of an observation of
Ferapontov \FerapontovL.

\Thm[Theorem]
{There is a local bijective correspondence between Lie minimal surfaces and
negatively curved projectively minimal surfaces modulo congruence.}

\vskip14pt{\bf Acknowledgements.}
The first author is grateful for the hospitality shown by Fukuoka
University, MPI Bonn and SFB288 Berlin during the preparation of this
paper.  Moreover, he thanks with pleasure V.~Cort\'es,
E.V.~Ferapontov, J.~Inoguchi and U.~Pinkall for instructive
conversations.


\References

\bye